\documentclass[11pt]{article}
\usepackage{latexsym}
\usepackage{amssymb}
\font\gorditas = msbm8
\def\bbb#1{\hbox {{\gordas #1}}}
\def\errita{\hbox{\gorditas R}}
\def\zetita{\hbox{\gorditas Z}}

\font\gordas = msbm10 at 12pt
\def\bbb#1{\hbox {{\gordas #1}}}
\def\erre{{\bbb R}}
\def\zet{{\bbb Z}}

\def\UNO{1\mkern-7mu1}

\newtheorem{theorem}{Theorem}[section]

\newtheorem{proposition}[theorem]{Proposition}

\newtheorem{remark}[theorem]{Remark}

\begin{document}

\begin{center}
{\bf\large Particle systems with quasi-homogeneous initial states and\\
their occupation time fluctuations$^*$}
 \vglue.5cm

\begin{tabular}{ccc}
Tomasz Bojdecki$^1$& Luis G. Gorostiza$^\dagger$&Anna  Talarczyk$^1$\\
tobojd@mimuw.edu.pl&lgorosti@math.cinvestav.mx&annatal@mimuw.edu.pl
\end{tabular}
\end{center}
\footnote{\kern-.6cm $^*$ Supported in part by CONACyT Grant 98998 (Mexico) and MNiSzW Grant N N201 397537 (Poland).\\
$^1$ Institute of Mathematics, University of Warsaw, ul. Banacha 2, 02-097 Warszawa, Poland.\\
$^\dagger$ Centro de Investigaci\'on y de Estudios Avanzados, A.P. 14-740, Mexico 07000 D.F., Mexico.}

\begin{abstract}
Occupation time fluctuation limits of particle systems in $\erre^d$ with independent motions (symmetric stable L\'evy process, with or without critical branching) have been studied assuming initial distributions given by Poisson random measures (homogeneous and some inhomogeneous cases). In this paper, with $d=1$ for simplicity, we extend previous results to a wide class of initial measures obeying a 
quasi-homogeneity property, which includes as special cases homogeneous Poisson measures and many deterministic measures (simple example: one atom at each point of $\zet$), by means of a new unified approach. In previous papers, in the homogeneous Poisson case, for the branching system in ``low'' dimensions, the limit was characterized by a long-range dependent Gaussian process called sub-fractional Brownian motion (sub-fBm), and this effect was attributed to the branching because it had appeared only in that case. An unexpected finding in this paper is that sub-fBm is more prevalent than previously thought. Namely, it is a natural ingredient of the limit process in the non-branching case (for ``low'' dimension), as well. On the other hand, fractional Brownian motion is not only related to systems in equilibrium (e.g., non-branching system with initial homogeneous Poisson measure), but it also appears here for a wider class of initial measures of quasi-homogeneous type.
\end{abstract}


\vskip.5cm

\newpage
\vglue.5cm

\section{Introduction}
\setcounter{equation}{0}

In a series of papers \cite{BGT0,BGT1,BGT2,BGT3,BGT4,BGT5,BGT6,BGT7} we studied particle systems in 
$\erre^d$ starting  from a configuration determined by a random point 
measure $\nu$, and independently moving according to a standard $\alpha$-stable L\'evy process ($0<\alpha\le 2$). In some models the particles additionally
 undergo critical branching. The evolution of the system is described by the empirical process $N=(N_t)_{t\geq 0}$, where 
$N_t(A)$ is the number of particles in the set $A\subset \erre^d$ at time $t$. 
The main object of   interest is the limit of the time-rescaled and normalized 
occupation time fluctuation process $X_T$ defined by

\begin{equation}
\label{eq:1.1}
X_T(t)=\frac{1}{F_T}\int_0^{Tt}(N_s-EN_s)ds, \quad t\geq 0,
\end{equation}
as $T\to\infty$ (i.e., as time is accelerated), where $F_T$ is a suitable deterministic norming. The process $X_T$ is signed measure-valued, but we regard it as a process with values in the space of tempered distributions 
${\cal S}'(\erre^d)$ for technical convenience, and also because in some cases the limit is genuinely ${\cal S}'(\erre^d)$-valued.
In all the cases considered in the abovementioned papers the initial measure 
$\nu$ was a Poisson field, homogeneous or not. This assumption permitted to investigate convergences conveniently with the help of the Laplace transform (due to infinite divisibility). The results always exhibited the same 
type of phase transition: for 
``low'' dimensions $d$ the limit process was the Lebesgue measure multiplied by a real long-range dependent process, whereas for ``high'' dimensions the limit was an ${\cal S}'(\erre^d)$-valued process with independent increments.

A natural question is what happens for  non-Poisson initial measures $\nu$. 
Mi\l o\'s \cite{Mi1,Mi2,Mi3} considered (critical) branching systems where 
$\nu$ was an  equilibrium measure (see \cite{GW}).
In that model the limits have a similar  dimension phase transition; moreover, for high dimensions they are the same as in the homogeneous Poisson case, while for low dimensions they are different.
 The conclusion was that in low dimensions the occupation time fluctuation process ``remembers'' the initial state of the system. Since the equilibrium states of the branching system are somewhat similar to homogeneous Poisson measures (they are infinitely divisible random point measures with uniform intensity; distributions of this kind are  called ``equilibrium distributions of Poisson type'' in \cite{LMW}), the Laplace transform method was also useful in 
\cite{Mi1,Mi2,Mi3}. 

The aim of the present paper is to investigate what happens with initial measures of other types, for example, some measures that are deterministic or almost deterministic. 
For simplicity we consider $d=1$ and assume that the motions are either without branching or with the simplest critical binary branching. 
In \cite{BGT1, BGT2} we proved for such motions,
 with general $d$, that if $\nu$ is a homogeneous Poisson measure, then  the following results hold (where $\lambda$ denotes Lebesgue measure and 
$K$ is a different constant in each case):

in the non-branching system:

if $d<\alpha$, then  $X_T$ converges in distribution (in $C([0,\tau],{\cal S}'(\erre^d))$ for any  $\tau>0$) to a process $K\lambda\xi$, where 
$\xi$ is a fractional Brownian motion;

if $d=\alpha$, then  the limit process is $K\lambda\beta$, where $\beta$ is a standard Brownian motion;

if $d>\alpha$,  then the limit is a time-homogeneous Wiener process in ${\cal S}'(\erre^d)$;

in the branching system:

if $\alpha<d<2\alpha$,  then the limit is $K\lambda\zeta$, where $\zeta$ is a sub-fractional Brownian motion (the case $d\le \alpha$ requires a slightly different treatment based on high-density models, see \cite{BGT6});

if $d=2\alpha$,  then the limit is $K\lambda\beta$;

if $d>2\alpha$,  then the limit is a time-homogeneous Wiener process in ${\cal S}'(\erre^d)$, different from the one  in the non-branching case.

In this paper we define a class $\cal M$ of initial measures $\nu$ which contains in particular homogeneous Poisson measures (which are ``completely random'' \cite{K}), and 
 quasi-homogeneous deterministic measures
 (e.g., the measure defined by one atom at each $j\in\zet$, which is ``completely deterministic''), and we develop a unified approach that permits to obtain limits of $X_T$ for all $\nu\in\cal M$. By a quasi-homogeneous deterministic point measure on $\erre$ we mean any measure defined by the following procedure: Given a positive integer $k$, in each interval $[j,j+1)$, $j\in\zet$, we fix $k$ points. For a general $\nu\in{\cal M}$, each interval $[j,j+1)$ contains $\theta_j$ points chosen at random,  and 
$\theta_j$, $j\in\zet$, are i.i.d.\ random variables (see Section 2 for a rigorous definition). The main feature of those measures is this form of quasi-homogeneity and independence on the family of  intervals $[j,j+1)$. 

For each $\nu\in{\cal M}$ we obtain the limit of the corresponding $X_T$ and in this way
 we recover  the results of \cite{BGT1, BGT2} for the homogeneous Poisson case (for $d=1$, but there is no doubt that the results for higher dimensions can be obtained analogously), and we also  derive
 limits for many other initial measures.
 It seems interesting that the  idea of the proofs in this general framework is simpler than that in our previous papers, and is based on the central limit theorem.  This is a significant change of methodology. However, some technical points in those papers are employed again here. An  argument using the non-linear equation associated with the occupation time of the branching system, which can be obtained by means of the Feynman-Kac theorem, again plays an important 
role, but now in a different way: it is a key step in moment estimates in order to apply the Lyapunov theorem in the branching case. The equilibrium measures for the branching system do not belong to ${\cal M}$ because the branching introduces spatial dependence.

Some of the results we obtain are unexpected. It turns out that the only case where new limits appear is the non-branching case with $(d=\,)1<\alpha$. They have the form $K\lambda\xi$, where $\xi$ is the sum of two independent processes, one of them is a sub-fractional Brownian motion (see (\ref{eq:2.3})), and the second one is a new (centered continuous with long range dependence) Gaussian process (see (\ref{eq:2.4})). The process $\xi$ depends on the initial measure $\nu$ only through $E\theta_0$ and ${\rm Var}\,\theta_0$.
%
In particular, for a  deterministic initial measure this process reduces  to a sub-fractional Brownian motion, and in the homogeneous Poisson case (as well as for any $\nu$ with $E\theta_0={\rm Var}\,
\theta_0$) it yields a fractional Brownian motion (see Theorem \ref{th:2.2}). 
This result seems surprising since in all earlier papers  sub-fractional Brownian motion was related only to  branching systems, and was consequently attributed to the branching, but now, in the present context, this process turns out to be more 
``natural'' than  fractional Brownian motion. On the other hand, fractional Brownian motion, which is typically related to systems in equilibrium (in particular the non-branching system with initial homogeneous Poisson measure), now appears also for a wider class of quasi-homogeneous initial measures, as noted above.

In all the remaining cases the limits are (up to constants)  the same, and with the same normings $F_T$, as those recalled above for homogeneous Poisson 
models.

The results show that within the class $\cal M$ the fluctuations caused by the branching are so large that $X_T$ ``forgets'' the randomness of the initial state of the system (it ``remembers'' $E \theta_0$ only). 
On the other hand, for low dimensions it does distinguish between $\nu\in{\cal M}$ and the equilibrium initial state (which is not in ${\cal M}$).
Another conclusion is that for high dimensions (which for $d=1$ amounts to small $\alpha$), the fluctuation process ``forgets''  the initial measure, as long as it is in some sense homogeneous (i.e., $\nu\in\cal M$), and this property holds for branching and non-branching systems; it is also preserved for branching systems in equilibrium.

In this paper we are interested mainly in identifying the limit processes, therefore we have not attempted to prove convergences in their strongest, functional form; in most cases we prove only convergence of finite-dimensional distributions. Presumably, convergence in distribution also holds in 
$C([0,\tau],{\cal S}'(\erre))$ for any $\tau>0$. As an example, we give one result of this type (Proposition \ref{prop:2.6}).

Some other papers related to occupation times of particle systems and  superprocesses are \cite{BZ, CG, DGW,   DR, DW, DF, FG, H,I,IL, Ta, Z}. For example, occupation time limits for branching random walks on the $d$-dimensional lattice are discussed in \cite{BZ}.

We have not found  results in the literature concerning occupation times for particle systems starting from a deterministic or quasi-deterministic point measure. 
Some kinds of quasi-homogeneity of initial configurations for systems of independent particles, different from those in this paper, appear in other contexts in 
\cite{S} and 
\cite{HJ} (see Remark \ref{rm:2.5} (f)). It may be that systems of independent particles with $\alpha$-stable motion and the initial conditions of \cite{S} lead to the same results as with initial homogeneous Poisson distribution.

The following notation is used in the paper.

${\cal S}(\erre)$: space of $C^\infty$ rapidly decreasing function on $\erre$.

${\cal S}'(\erre)$: space of tempered distributions (topological dual of ${\cal S}(\erre)$).

$\langle\,\, ,\,\,\rangle$: duality on ${\cal S}'(\erre)\times {\cal S}(
\erre)$.

$\Rightarrow_f$: weak convergence of finite-dimensional distributions of ${\cal S}'(\erre)$-valued processes.

$p_t(x)$: transition probability density of the standard symmetric $\alpha$-stable L\'evy process.

${\cal T}_t$:  semigroup determined by $p_t$, i.e., ${\cal T}_t\varphi=p_t*\varphi$.

$G$: potential operator determined by $p_t$ for $\alpha<1,$ i.e.,
\begin{equation}
\label{eq:1.2}
G\varphi(x)=\int^\infty_0{\cal T}_t\varphi(x)dt=C_\alpha\int_{\errita}\frac{\varphi(y)}{|x-y|^{1-\alpha}}dy,
\end{equation}
where
\begin{equation}
\label{eq:1.3}
C_\alpha=\frac{\Gamma(\frac{1-\alpha}{2})}{2^\alpha\pi^{1/2}\Gamma(\frac{\alpha}{2})}.
\end{equation}
Generic constants are written $C,C_i$, with possible dependencies in parenthesis.

In Section 2 we describe the particle system, formulate the results and discuss them. Section 3 contains the proofs.

\noindent
\section{Results}
\setcounter{section}{2}
\setcounter{equation}{0}

We start with detailed description of the particle system.

Let $\theta$ be a non-negative integer-valued random variable with distribution
\begin{equation}
\label{eq:2.1}
P(\theta=k)=p_k,\,\, k=0,1,2,\ldots,
\end{equation}
such that $E\theta^3<\infty$. This moment condition is a technical assumption satisfied by all cases of interest in this paper, but we suppose that finiteness of the second moment could be sufficient.

Let $\theta_j,j\in\zet$, be independent copies of $\theta$, and for each $j\in\zet$ and $k=1,2,\ldots$, let $\rho^j_k=(\rho^j_{k,1},\ldots,\rho^j_{k,k})$ be a random vector with values in $[j,j+1)^k$.
We assume that 
$(\theta_j,(\rho^j_k)_{k=1,2},\ldots),\, j\in \zet$, are independent.
These  objects determine a random point measure $\nu$ on $\erre$ in the following way: For each $j, \theta_j$ is the number of points in the interval $[j,j+1)$, and for each $k$, if $\theta_j=k$, the positions of those points are determined by $\rho^j_k$. In other words,
\begin{equation}
\label{eq:2.2}
\nu=\sum_{j\in\zetita}\sum^{\theta_j}_{n=1}
\delta_{\kappa_{j,n}},
\end{equation}
where $\kappa_{j,n}=\rho^j_{\theta_{j},n}$ and $\delta_a$ is the Dirac  measure at $a\in\erre$.
We denote by ${\cal M}$ the class of all such measures $\nu$.

\begin{remark}\label{rm:2.1}
{\rm (a) If $\theta\equiv k$ and, for each $j$, the $\rho^j_k$ are not random, then $\nu$ is a quasi-homogeneous deterministic measure mentioned in Introduction. The simplest example is $\nu=\sum_{j\in \zetita}\delta_j$.

(b) If $\theta$ is a standard Poisson random variable and, for each $j,\rho^j_{k,1},\ldots,\rho^j_{k,k}$ are independent, uniformly distributed on $[j,j+1)$, then $\nu$ given by (\ref{eq:2.2}) is the homogeneous Poisson point measure (with intensity measure $\lambda$).}
\end{remark}

Fix $\alpha\in(0,2]$ and $\nu\in{\cal M}$. Assume that at the initial time $t=0$ there is a collection of particles in $\erre$ with positions determined by $\nu$. As time evolves, these particles move independently according to the standard $\alpha$-stable L\'evy process. We consider systems either without branching, or with critical binary branching (i.e., $0$ or $2$ particles with probability $1/2$ each case) at rate $V$. For the corresponding empirical process, we define an ${\cal S}'(\erre)$-valued process $X_T$ by (\ref{eq:1.1}).

Before stating the first theorem we recall the definition of sub-fractional Brownian motion. A {\it  sub-fractional Brownian motion} with parameter $H\,(0<H<1)$ is a centered continuous Gaussian process $\zeta^H$ with covariance
\begin{equation}
\label{eq:2.3}
C^H(s,t)=E\zeta^H_s\zeta^H_t=s^{2H}+t ^{2H}-\frac{1}{2}((s+t)^{2H}+
|s-t|^{2H}),\,\, s, t\geq 0.
\end{equation}
See \cite{BGT0, T1} for properties of this process. It appears in \cite{DZ} in a different context, and it  has also been investigated in 
\cite{BD, RT, T2,T3,T4}.

We will need another centered Gaussian process $\vartheta^H$ with covariance
\begin{equation}
\label{eq:2.4}
Q^H(s,t)=\frac{1}{2}{\rm sgn}(2H-1)((s+t)^{2H}-s^{2H}-t^{2H}),\,\, s,t\geq 0, 
\end{equation}
$(0<H<1)$. Existence of this process for $H\neq 1/2$ follows from the formula
$$Q^H(s,t)=C\int^s_0\int^t_0(r+r')^{2H-2}dr'dr=C_1
\int^s_0\int^t_0\int_{\errita}e^{-r|x|^{1/(2-2H)}}e^{-r'|x|^{1/(2-2H)}}dxdr'dr,$$
which implies positive-definiteness of $Q^H$.

\begin{theorem}\label{th:2.2}
For the system without branching,

(a) if $1<\alpha$ and
\begin{equation}
\label{eq:2.5}
F_T=T^{1-1/2\alpha},
\end{equation}
then
\begin{equation}
\label{eq:2.6}
X_T\Rightarrow_fK_1\lambda(\sqrt{E\theta}\zeta^H+
\sqrt{{\rm Var}\,\theta}\vartheta^H),
\end{equation}
where $\zeta^H,\vartheta^H$ are independent, $H=1-1/2\alpha$, and
$$K_1=\left(
\frac{\Gamma(2-2H)}{2\pi\alpha H(2H-1)}\right)^{1/2};$$
(b) if $\alpha=1$ and
\begin{equation}
\label{eq:2.7}
F_T=\sqrt{T\log T},
\end{equation}
then
$$X_T\Rightarrow_fK_2\lambda\beta,$$
where $\beta$ is a standard Brownian motion in $\erre$, and
$$K_2=\sqrt{\frac{2}{\pi}E\theta};$$
(c) if $1>\alpha$ and
\begin{equation}
\label{eq:2.8}
F_T=\sqrt{T},
\end{equation}
then
\begin{equation}
\label{eq:2.9}
X_T\Rightarrow_f  X,
\end{equation}
where $X$ is an ${\cal S}'(\erre)$-valued homogeneous Wiener process with covariance
\begin{equation}
\label{eq:2.10}
E\langle X(t),\varphi\rangle\langle X(s),\psi\rangle=
2E\theta(s\wedge t)\int_{\errita}\varphi(x)G\psi(x)dx,\,\, \varphi,\psi\in{\cal S}(\erre).
\end{equation}
\end{theorem}

\begin{remark}\label{rm:2.3}
{\rm (a) Let 
\begin{equation}
\label{eq:2.11}
\xi^H=\sqrt{E\theta}\zeta^H+\sqrt{{\rm Var}\theta}\vartheta^H
\end{equation}
be the process in Theorem 2.2(a). If $E\theta={\rm Var}\theta$, in particular if $\nu$ is homogeneous Poisson (see Remark 2.1(b)), then $\xi^H$ is, up to a constant, a fractional Brownian motion with Hurst parameter $H$, i.e., it has covariance  $C(s^{2H}+t^{2H}-|s-t|^{2H})$. Thus we recover Theorem 2.1 of \cite{BGT1}. On the other hand, if $\theta$ is deterministic, then $\xi^H$ is a sub-fractional Brownian motion. Moreover, in general randomness of the $\rho$'s in the definition of $\nu\in{\cal M}$ (see (2.2)) does not play any role in the limit.

(b) The long time dependent behavior of Gaussian processes is usually characterized by the covariance of increments of the process on intervals separated by distance $\tau$, as $\tau\to\infty$. For the process $\vartheta^H$ that behavior is asymptotic decay like $\tau^{2H-2}$ (the same as for fractional Brownian motion), for sub-fBm it is  $\tau^{2H-3}$ (see \cite{BGT0}). So, the long time dependent behavior of the process $\xi^H$ in Theorem \ref{th:2.2}(a) is determined by $\vartheta^H$ in all cases where $\theta$ is random.
}

%
\end{remark}
\begin{theorem}\label{th:2.4}
For the branching system,

(a) if $1/2<\alpha<1$ and 
\begin{equation}
\label{eq:2.12}
F_T= T^{(3-1/\alpha)/2},
\end{equation}
then
$$X_T\Rightarrow_f K_3\lambda\zeta^H,$$
where $\zeta^H$ is a sub-fractional Brownian motion with parameter 
$H=(3-1/\alpha)/2$, and
$$K_3=\left(\frac{E\theta V\Gamma(2-2H)}{2\pi\alpha H(2H-1)}\right)^{1/2};$$

(b) if $\alpha=1/2$ and $F_T=\sqrt{T\log T}$, then
$$X_T\Rightarrow_fK_4\lambda\beta,$$
where $\beta$ is a standard Brownian motion, and
$$K_4=\sqrt{\frac{2V}{\pi}E\theta};$$
(c) if $\alpha<1/2$ and $F_T=\sqrt{T}$, then
$$X_T\Rightarrow_f X,$$
where $X$ is an ${\cal S}'(\erre)$-valued homogeneous Wiener process with covariance
\begin{equation}
\label{eq:2.13}
E\langle X(t),\varphi\rangle\langle X(s),\psi\rangle=E\theta(s\wedge t)\int_{\errita}\Big ( 2\varphi(x)G\psi(x)+V(G\varphi(x))(G\psi(x))\Big )dx,\,\, \varphi,\psi
\in{\cal S}(\erre).
\end{equation}
\end{theorem}

\begin{remark}\label{rm:2.5}{\rm 
(a) In the branching case the results are, up to the constant $E\theta$ in the limits, the same as in the homogeneous Poisson case (Theorems 2.2 in \cite{BGT1} and \cite{BGT2}).

(b) The condition $\alpha<1$ in part (a) of the last theorem corresponds to $\alpha<d$ in \cite{BGT1} and \cite{BGT2}.  In the homogeneous Poisson case for 
$d\leq \alpha$, we obtained  limits of the same form as for 
$\alpha<d<2\alpha$, by introducing high density, i.e., considering systems with initial intensity $H_T\lambda$,  
$H_T\to \infty$ sufficiently fast \cite{BGT6}; the high density counteracts the 
tendency to local extinction caused by the critical branching. The same procedure can be applied in the present case, yielding the limits for $1\leq\alpha\leq 2$, if the intervals $[j,j+1)$ are replaced by $[j/H_T, (j+1)/H_T)$.

(c) As in \cite{BGT1} and \cite{BGT2}, Theorems 2.2 and 2.4 can be extended to systems in $\erre^d$, where the intervals $[j,j+1)$ are replaced by cubes 
$[j,j+1)^d$.

(d) Comparing parts (a) of Theorems 2.2 and 2.4, we see that the branching weakens the influence of the initial configuration.

(e) The previous results show that sub-fractional Brownian motion is a ``natural'' process for our model. So far it had appeared only in the context of branching systems, but now we see that it is intrinsically related to the non-branching systems as well for a large class of initial conditions. Fractional Brownian motion occured before only in the case of systems in equilibrium, but now it also appears wherever $E\theta={\rm Var}\theta$.

(f) Theorems 2.2 and 2.4 can  also be extended to other models. For example, in \cite{HJ} a model is studied in a different context with independent $\alpha$-stable motions without branching and initial positions of particles 
$(j+\rho)_{j\in\zetita}$, where $\rho$ is a random variable uniformly distributed on $[0,1]$, independent of the motions. It is easy to see, by a standard conditioning argument (considering the characteristic function and conditioning on $\rho$), that for models with or without branching and with this initial configuration, the limits are the same as for  deterministic $\rho$, i.e., they are given by Theorems 2.2 and 2.4.  }
\end{remark}

We have formulated Theorems 2.2 and 2.4 with convergence of finite-dimensional distributions only, as we are mostly interested in the limit processes, but we have no doubt that  functional convergence also holds. As an example, let us consider the cases of large $\alpha$. For simplicity we assume that the initial configuration is such that $\theta\equiv 1$ and $\rho^j_{1,1}$ is uniformly distributed on $[j,j+1), j\in \zet$.

\begin{proposition}\label{prop:2.6}
For the model described above, the processes $X_T$  in Theorems 2.2(a) and 2.4(a) converge in law in $C([0,\tau], {\cal S}'(\erre))$ for any $\tau >0$. 
\end{proposition}

\section{Proofs}
\setcounter{section}{3}
\setcounter{equation}{0}

\subsection{Auxiliary facts related to the stable density}

\vglue.5cm
We will often use the  self-similarity property of the transition density $p_t$ of the standard $\alpha$-stable process in $\erre$:
\begin{equation}
\label{eq:3.1}
p_{at}(x)=a^{-1/\alpha}p_t(a^{-1/\alpha}x),\,\, x\in \erre, \,\, a>0.
\end{equation}
Recall that
\begin{equation}
\label{eq:3.2}
p_1(x)\leq\frac{C}{1+|x|^{1+\alpha}}.
\end{equation}
Since $p_t(\cdot)$ is decreasing on $\erre_+$ and symmetric, then by (\ref{eq:3.1}) we have
\begin{equation}
\label{eq:3.3}
p_t(x+y)\leq g_t(x):=
\left\{
\begin{array}{rll}
t^{-1/\alpha}p_1(0),&{\rm if}&|x|\leq 2\\
p_t(\frac{x}{2}),&{\rm if}&|x|>2
\end{array}\right. ,\quad x\in \erre,\,\,|y|\leq 1.
\end{equation}

Denote
\begin{equation}
\label{eq:3.4}
\phi_m(x)=\frac{1}{1+|x|^m},\,\, m>0.
\end{equation}
For $\varphi\in{\cal S}(\erre)$ we have $|\varphi(x)|\leq C(\varphi,m)\phi_m(x)$. This, and an obvious inequality,
\begin{equation}
\label{eq:3.5}
\frac{1}{1+|x+y|^m}\leq C(m)\frac{1+|y|^m}{1+|x|^m},\,\, m>0,
\end{equation}
imply
\begin{equation}
\label{eq:3.6}
{\cal T}_t|\varphi|(x+y)\leq C_m(a){\cal T}_t\phi_m(x),\,\, |y|\leq a, m>0,
\end{equation}
In the sequel we will use various versions of this estimate, e.g.,
$${\cal T}_t(|\varphi|{\cal T}_s|\psi|)(x+y)\leq C_m(a){\cal T}_t(\phi_m{\cal T}_s\phi_m)(x),\,\,\varphi,\psi\in {\cal S}(\erre).$$

We will also need the following estimate (\cite{I}, Lemma 5.3) for the potential operator $G$ (see (\ref{eq:1.2}),(\ref{eq:1.3})). If $\alpha<1, q>1$, and $\varphi$ is a measurable function on $\erre$ such that
$$|\varphi(x)|\leq C_1\frac{1}{1+|x|^q}, \,\, x\in \erre,$$
then
\begin{equation}
\label{eq:3.7}
|G\varphi(x)|\leq C_2\frac{1}{1+|x|^{1-\alpha}},\,\, x\in \erre.
\end{equation}

\subsection{Scheme of proofs}

The proofs of Theorems 2.2 and 2.4 are based on the central limit theorem and follow the scheme described presently.

Let $N^x$ denote the empirical process of the system (with or without branching) started from a single particle at $x$, and $N^{(j)}, j\in\zet$, be the empirical process for the particles which at time $t=0$ belong to $[j,j+1)$, i.e.,
\begin{equation}
\label{eq:3.8}
N^{(j)}=\sum^{\theta_j}_{n=1}N^{\kappa_{j,n}},
\end{equation}
according to the description at the beginning of Section 2 (see (2.2)). Note that $N^{(j)},j\in\zet$, are independent.

The process $X_T$ defined in (\ref{eq:1.1}) can be written as 
\begin{equation}
\label{eq:3.9}
X_T(t)=\sum_{j\in\zetita}\frac{1}{F_T}\int^{Tt}_0(N^{(j)}_s-EN^{(j)}_s)ds.
\end{equation}
The first step in our argument is to prove that for any 
$\varphi,\psi\in {\cal S}(\erre)$, and $s,t\geq 0$,
\begin{equation}
\label{eq:3.10}
\lim_{T\to\infty}E\langle X_T(t),\varphi\rangle\langle X_T(s),\psi\rangle=E\langle X(t),\varphi\rangle\langle X(s),\psi\rangle,
\end{equation}
where $X$ is the corresponding limit process. Without loss of generality we may assume that $\varphi,\psi\geq 0$.

Using (3.9) we have
\begin{eqnarray}
\label{eq:3.11}
\lefteqn{E\langle X_T(t),\varphi\rangle\langle X_T(s),\psi\rangle}\nonumber\\
&=&\sum_{j\in\zetita}\frac{1}{F^2_T}\int^{Tt}_0\int^{Ts}_0E\langle N^{(j)}_r,\varphi\rangle\langle N^{(j)}_{r'},\psi\rangle dr'dr\nonumber\\
&-&\sum_{j\in\zetita}\frac{1}{F^2_T}\int^{Tt}_0\int^{Ts}_0E\langle N^{(j)}_r,\varphi\rangle E\langle N^{(j)}_{r'},\psi\rangle dr'dr.
\end{eqnarray}
Using (\ref{eq:3.8}), (\ref{eq:2.1}) and the fact that 
$E\langle N^x_t,\varphi\rangle={\cal T}_t\varphi(x)$ in both non-branching and (critical) branching cases, and defining, for $x\in\erre$, $n\leq k$, random variables
\begin{equation}
\label{eq:3.12}
h_{k,n}(x)=\rho^{[x]}_{k,n}-x,
\end{equation}
where $[x]$ is the largest integer $\leq x$, we rewrite (\ref{eq:3.11}) as
\begin{eqnarray}
\label{eq:3.13}
\lefteqn{\kern-2cm E\langle X_T(t),\varphi\rangle\langle X_T(s),\psi\rangle=\sum^\infty_{k=0}p_k\sum^k_{n=1}I(T;k,n)+\sum^\infty_{k=0}p_k\sum^k_{n,m=1\atop_{n\neq m}}I\!\!I(T;k,n,m)}\nonumber\\
&-&\sum^\infty_{k=0}p_k\sum^k_{n=1}\sum^\infty_{\ell=0}p_\ell
\sum^\ell_{m=1}I\!\!I\!\!I(T;k,n;\ell,m),
\end{eqnarray}
where
\begin{eqnarray}
\label{eq:3.14}
I(T;k,n)&=&\frac{1}{F^2_T}\int^{Tt}_0\int^{Ts}_0\sum_{j\in\zetita}E\langle N^{\rho^j_{k,n}}_r,\varphi\rangle\langle N^{\rho^j_{k,n}}_{r'},\psi\rangle dr' dr\nonumber\\
&=&\frac{1}{F^2_T}\int^{Tt}_0\int^{Ts}_0\int_{\errita}E\langle N^{x+h_{k,n}(x)}_r,\varphi\rangle\langle N^{x+h_{k,n}(x)}_{r'},\psi\rangle dxdr'dr,\\
I\!\!I(T;k,n,m)&=&\frac{1}{F^2_T}\int^{Tt}_0\int^{Ts}_0\sum_{j\in\zetita}E\left({\cal T}_r\varphi(\rho^j_{k,n}){\cal T}_{r'}\psi(\rho^j_{k,m})\right)dr' dr\nonumber\\\label{eq:3.15}
&=&\frac{1}{F^2_T}\int^{Tt}_0\int^{Ts}_0\int_{\errita}E\left({\cal T}_r\varphi(x+h_{k,n}(x)){\cal T}_{r'}\psi(x+h_{k,m}(x))\right)dx dr' dr
\end{eqnarray}
(in the first equality for $I\!\!I$ we used independence of systems starting from different points),
\begin{eqnarray}
\label{eq:3.16}
I\!\!I\!\!I(T;k,n;\ell,m)&=&
\frac{1}{F^2_T}\int^{Tt}_0\int^{Ts}_0
\sum_{j\in\zetita}E{\cal T}_r\varphi(\rho^j_{k,n})
E{\cal T}_{r'}\psi(\rho^j_{\ell,m})dr' dr\nonumber\\
&=&\frac{1}{F^2_T}\int^{Tt}_0\!\!\int^{Ts}_0\int_{\errita}\!E{\cal T}_r\varphi(x+h_{k,n}(x))E{\cal T}_{r'}\psi(x+h_{\ell,m}(x))dx dr' dr.
\end{eqnarray}
Note that
\begin{equation}
\label{eq:3.17}
|h_{k,n}(x)|\leq 1,\,\, x\in\erre.
\end{equation}
In each case we will show convergence of $I, I\!\!I$ and $I\!\!I\!\!I$, thus proving (\ref{eq:3.10}). 
(It will be shown that $I, I\!\!I, I\!\!I\!\!I$ are bounded, so the passage to the limit in each sum in (\ref{eq:3.13}) is justified). 

Next, we show that 
$$\langle X(t),\varphi\rangle \Rightarrow \langle X(t),\varphi\rangle, \,\,\varphi\in{\cal S}(\erre), t\geq 0.$$
To this end, 
by (\ref{eq:3.9}) and (\ref{eq:3.10}) it suffices to prove that the Lyapunov condition
$$\lim_{T\to\infty}\sum_{j\in\zetita}\frac{1}{F^3_T}E\left|\int^{Tt}_0
(\langle N^{(j)}_r,\varphi\rangle-E\langle N^{(j)}_r,\varphi\rangle )
 dr\right|^3=0$$
is satisfied, and this property will follow if we show that
\begin{equation}
\label{eq:3.18}
\lim_{T\to\infty}\sum_{j\in\zetita}\frac{1}{F^3_T}E\left(\int^{Tt}_0\langle N^{(j)}_r,\varphi\rangle dr\right)^3=0,\,\, t\geq 0,\varphi\in {\cal S}(\erre),\varphi\geq 0.
\end{equation}

It is clear that convergence in law of linear combinations $\sum^m_{k=1}a_k\langle X_T(t_k),\varphi_k\rangle$ can be obtained analogously 
from (\ref{eq:3.10}) and (\ref{eq:3.18}), thus establishing the claimed convergence $X_T\Rightarrow_f X$.

In order to give (\ref{eq:3.18}) a more tractable form we use 
(\ref{eq:2.1}), (\ref{eq:3.8}), and the trivial inequality $(a_1+\ldots +a_k)^3\leq 3k^2(a^3_1+\ldots +a^3_k), a_1\ldots,a_k\geq 0$, obtaining
\begin{eqnarray*}
\lefteqn{\sum_{j\in\zetita}\frac{1}{F^3_T}
E\left(\int^{Tt}_0\langle N^{(j)}_r,\varphi\rangle dr\right)^3}\\
&=&\sum_{j\in\zetita}\frac{1}{F^3_T}\sum^\infty_{k=0}p_kE\left(\sum^k_{n=1}\int^{Tt}_0\langle N^{\rho^j_{k,n}}_r,\varphi\rangle dr\right)^3\\
&\leq&3\frac{1}{F^3_T}\sum^\infty_{k=0}p_kk^2\sum^k_{n=1}\sum_{j\in\zetita}E\left(\int^{Tt}_0\langle N^{\rho^j_{k,n}}_r,\varphi\rangle dr\right)^3\\
&\leq&3E\theta^3\sup_{n,k\in\zetita_+\atop_{n\leq k}}\frac{1}{F^3_T}\int_{\errita}E\left(\int^{Tt}_0\langle N^{x+h_{k,n}(x)}_r,\varphi\rangle dr\right)^3dx
\end{eqnarray*}
(see (\ref{eq:3.12})).
So, to prove (\ref{eq:3.18}) it suffices to show that
\begin{equation}
\label{eq:3.19}
\lim_{T\to\infty}\sup_{n,k\in\zetita_+\atop_{n\leq k}}\frac{1}{F^3_T}\int_{\errita}E\left(\int^{Tt}_0\langle N^{x+h_{k,n}(x)}_r,\varphi\rangle dr\right)^3dx=0,\,\, t\geq 0,\varphi\in{\cal S}(\erre),\varphi\geq 0.
\end{equation}

Summarizing, to obtain Theorems 2.2 and 2.4 we prove convergences of 
(\ref{eq:3.14})-(\ref{eq:3.16}) and (\ref{eq:3.19}). 
In each case these proofs require some non-trivial work.

\subsection{Proof of Theorem 2.2(a)}

Following the scheme we show that
\begin{eqnarray}
\label{eq:3.20}
\lim_{T\to\infty}I(T;k,n)&=&p_1(0)\int^t_0\int^s_0|r-r'|^{-1/\alpha}dr'dr\int_{\errita}\varphi(x)dx\int_{\errita}\psi(x)dx,\\
\lim_{T\to\infty}I\!\!I(T;k,n,m)&=&\lim_{T\to\infty}I\!\!I\!\!I(T;k,n;\ell,m)\nonumber\\
\label{eq:3.21}
&=&p_1(0)\int^t_0\int^s_0(r+r')^{-1/\alpha}dr'dr\int_{\errita}\varphi(x)dx\int_{\errita}\psi(x)dx.
\end{eqnarray}
(see (\ref{eq:3.14})-(\ref{eq:3.16})). It is easy to see that by 
(\ref{eq:3.13}), (\ref{eq:2.3}) and (\ref{eq:2.4}), this yields (\ref{eq:3.10}).

Let $\eta$ denote the standard $\alpha$-stable L\'evy process in $\erre$. As we consider the model without branching, we have, for $r>r'$,
\begin{eqnarray}
\lefteqn{E\langle N^{x+h_{k,n}(x)}_r,\varphi\rangle \langle N^{x+h_{k,n}(x)}_{r'},\psi\rangle}\nonumber\\
&=&E\varphi(x+h_{k,n}(x)+\eta_r)\psi(x+h_{k,n}(x)+\eta_{r'})\nonumber\\
\label{eq:3.22}
&=&{\cal T}_{r'}(\psi{\cal T}_{r-r'}\varphi)(x+h_{k,n}(x)).
\end{eqnarray}
Putting this into (\ref{eq:3.14}) and omitting the subscripts $k,n$, we obtain 
\begin{equation}
\label{eq:3.23}
I(T)=I_1(T)+I_2(T),
\end{equation}
where
\begin{eqnarray}
\label{eq:3.24}
I_1(T)&=&\frac{1}{F^2_T}\int^{Tt}_0\!\!\!\int^{Ts}_0\UNO_{\{r>r'\}}\!\!\int_{\errita^3}p_{r'}(x+h(x)-y)\psi(y)p_{r-r'}(y-z)\varphi(z)dzdydxdr'dr,
\kern1cm \\
\label{eq:3.25}
I_2(T)&=&\frac{1}{F^2_T}\int^{Tt}_0\!\!\!\int^{Ts}_0\UNO_{\{r\leq r'\}}\!\!\int_{\errita^3}p_{r}(x+h(x)-y)\varphi(y)p_{r'-r}(y-z)\psi(z)dzdydxdr'dr.
\end{eqnarray}

In $I_1$ we substitute $\widetilde{r}=r/T,\widetilde{r}'=r'/T$, use
(\ref{eq:2.5}) and (\ref{eq:3.1}), and then we substitute $\widetilde{x}=T^{-1/\alpha}(x-y)$, arriving at
\begin{eqnarray*}
I_1(T)
&=&\int^t_0\int^s_0\UNO_{\{r>r'\}}\int_{\errita^3}p_{r'}(x+T^{-1/\alpha}h(T^{1/\alpha}x+y))\psi(y)\\
&&\cdot p_{r-r'}(T^{-1/\alpha}(y-z))\varphi(z)dzdydxdr'dr.
\end{eqnarray*}
By (\ref{eq:3.17}) and (\ref{eq:3.1}), the expression under the integrals converges pointwise, as $T\to\infty$, to
$$p_{r'}(x)p_{r-r'}(0)\psi(y)\varphi(z)=p_{r'}(x)
|r-r'|^{-1/\alpha}p_1(0)\psi(y)\varphi(z),$$
and by (\ref{eq:3.3}), for $T>1$, it is bounded by 
$g_{r'}(x)|r-r'|^{-1/\alpha}p_1(0)\psi(y)\varphi(z)$, which is integrable over $[0,t]\times [0,s]\times \erre^3$, since $\alpha>1$. $I_2(T)$ can be treated analogously, hence by (\ref{eq:3.23}) we obtain (\ref{eq:3.20}).

Next we take $I\!\!I$. In (\ref{eq:3.15}) we substitute 
$\widetilde{r}=r/T,\widetilde{r}'=r'/T$, and by (\ref{eq:2.5}) we have
\begin{eqnarray*}
\lefteqn{I\!\!I(T;k,n,m)}\\
&=&T^{1/\alpha}E\int^t_0\int^s_0\int_{\errita^3}p_{Tr}(x+h_{k,n}(x)-y)p_{Tr'}(x+h_{k,m}(x)-z)\varphi(y)\psi(z)dxdydzdr'dr.
\end{eqnarray*}
We use (\ref{eq:3.1}) and substitute 
$\widetilde{x}=T^{-1/\alpha}(x-z)$, obtaining
\begin{eqnarray*}
I\!\!I(T;k,n,m)&=&E\int^t_0\int^s_0\int_{\errita^3}p_r(x+T^{-1/\alpha}(z-y+h_{k,n}(T^{1/\alpha}x+z))\\
&&\hskip 1.7cm\cdot p_{r'}(x+T^{-1/\alpha}h_{k,m}(T^{1/\alpha}x+z))\varphi(y)\psi(z)dxdydzdr'dr.
\end{eqnarray*}
The integrand converges pointwise to $p_r(x)p_{r'}(x)\varphi(y)\psi(z)$, and 
(\ref{eq:3.17}) implies that for $T>1$ , it is bounded by 
$r^{-1/\alpha}p_1(0)g_{r'}(x)\varphi(y)\psi(z)$ (see (\ref{eq:3.3})), which is integrable. As $\int_{\errita}p_r(x)p_{r'}(x)dx=(r+r')^{-1/\alpha}p_1(0)$, 
we obtain the limit (\ref{eq:3.21}) for $I\!\!I$.

Note that in this argument the only property of $h_{k,n}$ we have used is 
(\ref{eq:3.17}), therefore it is immediately seen that the limit of 
$I\!\!I\!\!I$ can be obtained in the same way (see (\ref{eq:3.16})).  This completes the proof of 
(\ref{eq:3.10}).

As explained in the scheme, to finish the proof it suffices to show (\ref{eq:3.19}). 

The expression under $\lim_{T\to\infty}\sup_{n,k}$ in (\ref{eq:3.19}), similarly as in (\ref{eq:3.22}), can be written as

\begin{eqnarray}
\lefteqn{\frac{3!}{F^3_T}
\int_{\errita}
\int^{Tt}_0\int^{Tt}_r\int^{Tt}_{r'}
{\cal T}_r(\varphi{\cal T}_{r'-r}(\varphi{\cal T}_{r''-r'}
\varphi))(x+h_{k,n}(x))dr''dr'drdx}\nonumber\\
&\leq&\frac{6}{F^3_T}
\int_{\errita}\int^{Tt}_0\int^{Tt}_0\int^{Tt}_0
{\cal T}_r(\varphi{\cal T}_{r'}(\varphi{\cal T}_{r''}
\varphi))(x+h_{k,n}(x))dr''dr'drdx\nonumber\\
\label{eq:3.26}
&\leq&C J(T),
\end{eqnarray}
where
\begin{equation}
\label{eq:3.27}
J(T)=\frac{1}{F^3_T}\int_{\errita}\int^{Tt}_0\int^{Tt}_0\int^{Tt}_0
{\cal T}_r(\phi_2{\cal T}_{r'}(\phi_2{\cal T}_{r''}\phi_2))(x)dr''dr'drdx,
\end{equation}
and $\phi_2$ is given by (\ref{eq:3.4}). The last inequality in (\ref{eq:3.26}) is of the same type as (\ref{eq:3.6}), and can be obtained by an analogous argument using (\ref{eq:3.5}) and (\ref{eq:3.17}). Hence, for (\ref{eq:3.19}) it is enough to show that
\begin{equation}
\label{eq:3.28}
\lim_{T\to\infty}J(T)=0.
\end{equation}
Note that in this argument we have not used the assumption on $\alpha$.

After obvious substitutions, using (\ref{eq:2.5}) and the invariance of Lebesgue measure for ${\cal T}_t$, we have
$$J(T)=T^{3/2\alpha}\int^t_0\int^t_0\int^t_0\int_{\errita^3}\phi_2(x)
p_{Tr'}(x-y)\phi_2(y)p_{Tr''}(y-z)\phi_2(z)dzdydxdr''dr'dr.$$
By the self-similarity property (\ref{eq:3.1}),
$$J(T)\leq T^{3/2\alpha-2/\alpha} p^3_1(0)t
\left(\int^t_0s^{-1/\alpha}ds\right)^2
\left(\int_{\errita}\phi_2(x)dx\right)^3\to 0 \,\,{\rm as}\,\,T\to\infty,$$
since $1<\alpha. \hfill\Box$

\subsection{Proof of Theorem 2.2 (b),(c)}

According to the scheme, we prove (\ref{eq:3.10}). Fix $s\leq t$. (\ref{eq:3.14}) can be written as
\begin{equation}
\label{eq:3.29}
I(T;k,n)=I_1(T;k,n)+I_2(T;k,n)+I_3(T;k,n),
\end{equation}
where
\begin{eqnarray}
\label{eq:3.30}
I_1(T;k,n)&=&\frac{1}{F^2_T}E\int^{Ts}_0\int^{Ts}_r\int_{\errita}{\cal T}_r(\varphi{\cal T}_{r'-r}\psi)(x+h_{k,n}(x))dxdr'dr,\\
\label{eq:3.31}
I_2(T;k,n)&=&\frac{1}{F^2_T}E\int^{Ts}_0\int^r_0\int_{\errita}{\cal T}_{r'}(\psi{\cal T}_{r-r'}\varphi)(x+h_{k,n}(x))dxdr'dr,\\
\label{eq:3.32}
I_3(T;k,n)&=&\frac{1}{F^2_T}E\int^{Tt}_{Ts}\int^{Ts}_0\int_{\errita}\ldots dxdr'dr,
\end{eqnarray}
where $\ldots$ is the same integrand as in (\ref{eq:3.31}).

To compute the limit of $I_1$ we fix an arbitrary $0<\varepsilon<1$, substitute $\tilde{r}=r/T$, and we write
\begin{equation}
\label{eq:3.33}
I_1(T;k,n)=I'_1(T;k,n)+I''_1(T;k,n),
\end{equation}
\begin{eqnarray}
\label{eq:3.34}
I'_1(T;k,n)&=&\frac{T}{F^2_T}E\int^{\varepsilon s}_0\int^{Ts}_{Tr}\int_{\errita}{\cal T}_{Tr}(\varphi{\cal T}_{r'-Tr}\psi)(x+h_{k,n}(x))dxdr'dr,\\
\label{eq:3.35}
I''_1(T;k,n)&=&\frac{T}{F^2_T}E\int^s_{\varepsilon s}\int^{Ts}_{Tr}\int_{\errita} \ldots dxdr'dr.
\end{eqnarray}

Consider case (c): $\alpha<1$.
Applying a version of (\ref{eq:3.6}) we have
\begin{eqnarray}
I'_1(T;k,n)&\leq&\frac{T}{F^2_T}C\int^{\varepsilon s}_0\int^{Ts}_0\int_{\errita}{\cal T}_{Tr}(\phi_2{\cal T}_{r'}\phi_2)(x)dxdr'dr\nonumber\\
\label{eq:3.36}
&=&\frac{T}{F^2_T}C\varepsilon s\int^{Ts}_0\int_{\errita}\phi_2(x){\cal T}_{r'}\phi_2(x)dxdr'\\
\label{eq:3.37}
&\leq&C_1\varepsilon\int_{\errita}\phi_2(x)G\phi_2(x)dx=C_2\varepsilon,
\end{eqnarray}
by (\ref{eq:3.7}) and (\ref{eq:2.8}). Next,
\begin{eqnarray}
\kern-2cm I''_1(T;k,n)&=&\frac{T}{F^2_T}E\int^s_{\varepsilon s}\int^{T(s-r)}_0\int_{\errita^2}p_{Tr}(x+h_{k,n}(x)-y)\varphi(y){\cal T}_{r'}\psi(y)dydxdr'dr
\label{eq:3.38}
\\
&=&E\int^s_{\varepsilon s}\int_{\errita^2}p_r(x+T^{-1/\alpha}h_{k,n}(T^{1/\alpha}x+y))\varphi(y)\int^{T(s-r)}_0{\cal T}_{r'}\psi(y)dr'dydxdr,\nonumber
\end{eqnarray}
where we have applied (\ref{eq:3.1}) and the substitution 
$\tilde{x}=T^{-1/\alpha}(x-y)$.

It is now easy to see that
\begin{equation}
\label{eq:3.39}
\lim_{T\to\infty}I''_1(T;k,n)=(s-\varepsilon s)\int_{\errita}\varphi(y)G\psi(y)dy,
\end{equation}
since the passage of the limit under the integrals is justified by (\ref{eq:3.3}). Note that the function $g_r(x)$ is not integrable in $r$ in a neighborhood of $0$ for $\alpha\leq 1$; that is why we had to consider the interval $[\varepsilon s,s]$.

From (\ref{eq:3.33}), (\ref{eq:3.37}) and (\ref{eq:3.39}), we infer that 
\begin{equation}
\label{eq:3.40}
\lim_{T\to\infty}I_1(T;k,n)=s\int_{\errita}\varphi(y)G\psi(y)dy.
\end{equation}
Note that (\ref{eq:3.31}) is essentially the same as (\ref{eq:3.30}), only the roles of $\varphi$ and $\psi$ are interchanged. Therefore, by symmetry of $G$, we have
\begin{equation}
\label{eq:3.41}
\lim_{T\to\infty}I_2(T;k,n)=s\int_{\errita}\varphi(y)G\psi(y)dy.
\end{equation}

Passing to $I_3$ defined by 
(3.32), we first estimate it similarly as $I_1$ 
(see(\ref{eq:3.36})), then we change the order of integration $drdr'$, and substitute $\tilde{r'}=r'/T$ and $\tilde{r}=r-T\tilde{r'}$, obtaining
$$I_3(T;k,n)\leq C\int^s_0\int^\infty_{T(s-r')}\int_{\errita}\phi_2(x){\cal T}_r\phi_2(x)dxdrdr'.$$
Now, it is clear that
$$\lim_{T\to\infty}I_3(T;k,n)=0.$$
hence, by (\ref{eq:3.29}), (\ref{eq:3.40}) and (\ref{eq:3.41}), we have
\begin{equation}
\label{eq:3.42}
\lim_{T\to\infty}I(T;k,n)=2(s\wedge t)\int_{\errita}\varphi(x)G\psi(x)dx.
\end{equation}

Next, by (\ref{eq:3.15}), (\ref{eq:3.16}), using (\ref{eq:3.17}) and 
(\ref{eq:3.6}), we obtain
\begin{eqnarray}
I\!\!I(T;k,n,m)&+&I\!\!I\!\!I(T;k,n;l,m)\leq\frac{C}{F^2_T}\int^{Tt}_0\int^{Ts}_0\int_{\errita}{\cal T}_r\phi_2(x){\cal T}_{r'}\phi_2(x)dxdr'dr\nonumber\\
\label{eq:3.43}
&=&\frac{C}{F^2_T}\int^{Tt}_0\int^{Ts}_0\int_{\errita}\phi_2(x)
{\cal T}_{r+r'}\phi_2(x)dxdr'dr\\
\label{eq:3.44}
&\leq&C\int^t_0\int^\infty_{Tr}\int_{\errita}\phi_2(x){\cal T}_{r'}\phi_2(x)dxdr'dr\to 0,
\end{eqnarray}
since, by (\ref{eq:3.7}), $G\phi_2$ is bounded. This together with (\ref{eq:3.42}) and (\ref{eq:3.13}) proves (\ref{eq:3.10}) for $\alpha <1$.

Now consider  case (b): $\alpha=1$. We also have 
(\ref{eq:3.36}), hence we get
$$I'_1(T;k,n)\leq C\varepsilon,$$
because for such $\alpha$,
\begin{equation}
\label{eq:3.45}
\sup_{T>2}\sup_{x\in\errita}\frac{1}{\log T}\int^T_0{\cal T}_r\phi_2(x)dr<\infty,
\end{equation}
see (3.46) in \cite{BGT5}. It is easy to see that the limit of $I''_1(T;k,n)$ given by (\ref{eq:3.38}) is the same as the limit of
$$\frac{1}{\log T}E\int^{s-1/T}_{\varepsilon s}\int^{T(s-r)}_1\int_{\errita^3}p_{Tr}(x+h_{k,n}(x)-y)\varphi(y)p_{r'}(y-z)\psi(z)dzdydxdr'dr.$$
Using (\ref{eq:3.1}) twice and substituting 
$\tilde{r'}=\log r'/\log T$ and $\tilde{x}=T^{-1}(x-y)$, this expression is equal to
\begin{eqnarray*}
\lefteqn{\kern-4cm E\int^{s-1/T}_{\varepsilon s}\int^{\log T(s-r)/\log T}_0\int_{\errita^3}p_r(x+T^{-1}h(Tx+y))\varphi(y)p_1(T^{-r'}(y-z))\psi(z)dzdydxdr'dr}\\
&\longrightarrow&s(1-\varepsilon)p_1(0)\int_{\errita}
\varphi(y)dy\int_{\errita}\psi(z)dz,
\end{eqnarray*}
by (\ref{eq:3.3}). Thus, we have shown that 
\begin{equation}
\label{eq:3.46}
\lim_{T\to\infty}I_1(T;k,n)=(s\wedge t)p_1(0)\int_{\errita}\varphi(x)dx\int_{\errita}\psi(y)dy=\lim_{t\to\infty}I_2(T;k,n),
\end{equation}
where the last equality follows by symmetry (see (\ref{eq:3.31})).

To estimate $I_3$ given by (3.32), we again use an inequality of the type (\ref{eq:3.6}), obtaining
\begin{eqnarray*}
I_3(T;k,n)&\leq&\frac{C}{T\log T}\int^{Tt}_{Ts}\int^{Ts}_0\int_{\errita^3}p_{r'}(x-y)\phi_2(y)p_{r-r'}(y-z)\phi_2(z)dzdydxdr'dr\\
&\leq&\frac{C}{T\log T}p_1(0)\left(\int_{\errita}\phi_2(y)dy\right)^2\int^{Tt}_{Ts}\int^{Ts}_0(r-r')^{-1}dr'dr\to 0\,\,{\rm as}\,\, T\to\infty.
\end{eqnarray*}
Similarly, by (\ref{eq:3.43}) we have
\begin{eqnarray*}
I\!\!I(T;k,n,m)&+&I\!\!I\!\!I(T;k,n;\ell,m)\\
&\leq&\frac{C}{\log T}p_1(0)\int^t_0\int^s_0(r+r')^{-1}dr'dr
\left(\int_{\errita}\phi_2(x)dx\right)^2\\
&&\longrightarrow 0\,\,\quad{\rm as}\,\, T\to\infty.
\end{eqnarray*}
This, together with (\ref{eq:3.46}), (\ref{eq:3.29}) and (\ref{eq:3.13}) prove 
(\ref{eq:3.10}).

It remains to show (\ref{eq:3.19}). To this end, by (\ref{eq:3.26}) and (\ref{eq:3.27}), it suffices to observe that (\ref{eq:3.28}) holds. Indeed, for $\alpha <1$ we use boundedness of $G\phi_2$, and for $\alpha=1$ we employ (\ref{eq:3.45}). The proof of the theorem is complete. $\hfill\Box$

\subsection{Proof of Theorem 2.4}

We recall first the following formula for the second moments of critical binary branching systems with branching rate V:
\begin{eqnarray}
\lefteqn{E\langle N^x_r,\varphi\rangle \langle N^x_{r'},\psi\rangle}\nonumber\\
\label{eq:3.47}
&=&{\cal T}_r(\varphi{\cal T}_{r'-r}\psi)+V\int^r_0{\cal T}_u(({\cal T}_{r-u}\varphi)({\cal T}_{r'-u}\psi))du,\,\, r\leq r', x\in \erre,
\end{eqnarray}
which is obtained using, e.g. Lemma 3.1 in  \cite{K1}, and the Markov property.

Fix $s\leq t$. By (\ref{eq:3.47}), expression (\ref{eq:3.14}) can be written as
\begin{equation}
\label{eq:3.48}
I(T;k,n)=VI_1(T;k,n)+I_2(T;k,n),
\end{equation}
where
\begin{eqnarray}
\label{eq:3.49}
I_1(T;k,n)&=&\frac{1}{F^2_T}E\int^{Ts}_0\int^{Tt}_u\int^{Ts}_u\int_{\errita}{\cal T}_u(({\cal T}_{r-u}\varphi)({\cal T}_{r'-u}\psi))(x+h_{k,n}(x))dxdr'drdu,\\
I_2(T;k,n)&=&\frac{1}{F^2_T}E\int^{Tt}_0\int^{Ts}_0\int_{\errita}\,\,({\cal T}_r(\varphi{\cal T}_{r'-r}\psi)(x+h_{k,n}(x))\UNO_{\{r\leq r'\}}\nonumber\\
&&\qquad\qquad\qquad\qquad+{\cal T}_{r'}(\psi{\cal T}_{r-r'}\varphi)(x+h_{k,n}(x))\UNO_{\{r'<r\}})dxdr'dr.\nonumber
\end{eqnarray}
Comparing $I_2$ to (\ref{eq:3.29})-(\ref{eq:3.32}), we see that it has the same form as $I(T;k,n)$ in the non-branching case with $\alpha<1$, only the norming $F_T$ may be different. It is now immediately seen that by (\ref{eq:3.42}) we have
\begin{equation}
\label{eq:3.50}
\lim_{T\to\infty}I_2(T;k,n)=\left\{\begin{array}{cll}
0,&{\rm if}&\frac{1}{2}\leq\alpha<1,\\
2(s\wedge t)\int_{\errita}\varphi(x)G\psi(x)dx,&{\rm if}&0<\alpha<\frac{1}{2}.
\end{array}\right.
\end{equation}

Similarly, using the result for the non-branching system (see (\ref{eq:3.44})), we obtain
\begin{equation}
\label{eq:3.51}
\lim_{T\to\infty}I\!\!I(T;k,n,m)=\lim_{T\to\infty}I\!\!I\!\!I(T;k,n;\ell,m)=0.
\end{equation}

Hence, to prove (\ref{eq:3.10}) it suffices to calculate the limit of 
$I_1(T;k,n)$ given by (\ref{eq:3.49}). 

Consider the case ${1}/{2}<\alpha<1$.  Making the substitutions 
$\tilde{r}=(r-u)/T,\tilde{r'}=(r'-u)/T, \tilde{u}=u/T$, using (\ref{eq:3.1}) 
and (\ref{eq:2.12}), and then putting $\tilde{x'}=T^{-1/\alpha}x, \tilde{y}=T^{-1/\alpha}y$, we have 
\begin{eqnarray}
\lefteqn{I_1(T;k,n)}\nonumber\\
&=&T^{-2/\alpha}E\int^s_0\int^{t-u}_0\int^{s-u}_0\int_{\errita^4}p_u(T^{-1/\alpha}(x+h_{k,n}(x)-y))p_r(T^{-1/\alpha}(y-z))\nonumber\\
\label{eq:3.52}
&&\qquad\qquad \cdot p_{r'}(T^{-1/\alpha}(y-w))\varphi(z)\psi(w)dwdzdydxdr'drdu\\
&=&E\int^s_0\int^{t-u}_0\int^{s-u}_0R_{T,k,n}(u,r,r')dr'drdu,\nonumber
\end{eqnarray}
where
\begin{eqnarray}
\label{eq:3.53}
R_{T,k,n}(u,r,r')
&=&\int_{\errita^4}p_u(x-y+T^{-1/\alpha}h_{k,n}(T^{1/\alpha}x))p_r(y-T^{-1/\alpha}z)\\
&&\cdot p_{r'}(y-T^{-1/\alpha}w)\varphi(z)\psi(w)dwdzdydx.\nonumber
\end{eqnarray}
Using (\ref{eq:3.2}), $\varphi\leq C\phi_m$, (\ref{eq:3.5}) and (\ref{eq:3.1}), it is easy to see that 
\begin{eqnarray*}
\lim_{T\to\infty}R_{T,k,n}(u,r,r')&=&\int_{\errita^4}p_u(x-y)p_r(y)p_{r'}(y)\varphi(z)\psi(w)dwdzdydx\\
&=&(r+r')^{-1/\alpha}p_1(0)\int_{\errita}\varphi(z)dz\int_{\errita}\psi(w)dw.
\end{eqnarray*}
By 
(\ref{eq:3.52}), (\ref{eq:3.51}), (\ref{eq:3.48}), (\ref{eq:3.50}), (\ref{eq:3.13}) and (\ref{eq:2.3}), it is clear that in order to obtain (\ref{eq:3.10}) it remains to justify the passage to the limit under the integral in (\ref{eq:3.52}). To this end, in (\ref{eq:3.53}) we substitute 
$\tilde{y}=y-T^{-1/\alpha}h_{k,n}(T^{1/\alpha}x),\tilde{z}=z-h_{k,n}
(T^{1/\alpha}x), \tilde{w}=w-h_{k,n}(T^{1/\alpha}x)$, use 
$\varphi(z+h_{k,n}(T^{1/\alpha}x))\leq C\phi_2(z)$ (see (\ref{eq:3.17}) and (\ref{eq:3.5})), and the same for $\psi$, obtaining
\begin{eqnarray*}
R_{T,k,n}(u,r,r')&\leq&C\int_{\errita^2}p_{r+r'}(T^{-1/\alpha}(z-w))\phi_2(z)\phi_2(w)dwdz\\
&\leq&C(r+r')^{-1/\alpha}p_1(0)\left(\int_{\errita}\phi_2(z)dz\right)^2.
\end{eqnarray*}
The last expression is integrable on the set of integration in (\ref{eq:3.52}), since $\alpha>1/2$.

Now consider  the case $0<\alpha<{1}/{2}$. After obvious substitutions, and using 
$F^2_T=T, I_1$ given by (\ref{eq:3.49}) can be written as
\begin{eqnarray}
I_1(T;k,n)
&=&E\int^s_0\int^{T(r-u)}_0\int^{T(s-u)}_0\int_{\errita}{\cal T}_{Tu}
(({\cal T}_r\varphi)({\cal T}_{r'}\psi))(x+h_{k,n}(x))dxdr'drdu\nonumber\\[.3cm]\label{eq:3.54}
&=&I'_1(T;k,n)+I''_1(T;k,n),
\end{eqnarray}
where, for any $0<\varepsilon<1$,
\begin{eqnarray}
\label{eq:3.55}
I'_1(T;k,n)&=&E\int^{\varepsilon s}_0 \ldots,\\
\label{eq:3.56}
I''_1(T;k,n)&=&E\int^s_{\varepsilon s}\ldots ,
\end{eqnarray}
(cf. (\ref{eq:3.33})-(\ref{eq:3.35})). We have
\begin{eqnarray*}
I'_1(T;k,n)&\leq&E\int^{\varepsilon s}_0
\int_{\errita^2}p_{Tu}(y)G\varphi(x+h_{k,n}(x)-y)G\psi(x+h_{k,n}(x)-y)dydxdu\\
&\leq&C\int^{\varepsilon s}_0
\int_{\errita^2}p_{Tu}(y)\frac{1}{1+|x-y|^{2(1-\alpha)}} dydxdu,
\end{eqnarray*}
by (\ref{eq:3.7}), (\ref{eq:3.5}) and (\ref{eq:3.17}), hence
\begin{equation}
\label{eq:3.57}
I'_1(T;k,n)\leq C_1\varepsilon,
\end{equation}
since $\alpha<1/2$.

Similarly as in (\ref{eq:3.38}) we obtain
$$\lim_{T\to\infty}I''_1(T;k,n)=(s-\varepsilon s)\int_{\errita}G\varphi(y)G\psi(y)dy,$$
hence, by (\ref{eq:3.57}) and (\ref{eq:3.54}),
$$\lim_{T\to\infty}I_1(T;k,n)=(s\wedge t)
\int_{\errita}G\varphi(y)G\psi(y)dy.
$$
This, together with (\ref{eq:3.48}), (\ref{eq:3.50}), (\ref{eq:3.51}), 
(\ref{eq:3.13}) and (\ref{eq:2.13}), proves (\ref{eq:3.10}).

It remains to show (\ref{eq:3.10}) in the case $\alpha=1/2$, which amounts to proving that
$$\lim_{T\to\infty}I_1(T;k,n)=(s\wedge t)\frac{2}{\pi}\int_{\errita}\varphi(x)dx\int_{\errita}\psi(x)dx.$$

We decompose $I_1$ as in (\ref{eq:3.54})-(\ref{eq:3.56}). The integral 
$(T/F^2_T)E\int^{\varepsilon s}_0\ldots$ is easy to estimate. To obtain the limit of $(T/F^2_T)E\int^s_{\varepsilon s}\ldots$ we combine the argument on pp. 1354-1356 of \cite{BGT7} and the method for getting rid of $h_{k,n}(x)$, which we have used several times above. We skip the cumbersome details.

We have proved (\ref{eq:3.10}) in all the cases. According to the scheme, to finish the proof it remains to show (\ref{eq:3.19}). To this and we define the function 
$$v_\theta(x,t)=1-E{\rm exp}
\left\{-\theta\int^t_0\langle N^x_r,\varphi\rangle dr\right\},
\,\, \theta\geq 0, x\in\erre, t\geq 0.$$
It is known that the Feynman-Kac formula implies that $v_\theta$ satisfies the non-linear equation
$$v_\theta(x,t)=\int^t_0{\cal T}_{t-u}
\left(\theta\varphi(\cdot)(1-v_\theta(\cdot,u))-\frac{V}{2}v_\theta^2
(\cdot,u)\right)(x)du$$
(see e.g. \cite{GLM}, or  the space-time approach used in 
 \cite{BGT1,BGT2}).
Hence, by a similar argument as in (3.45)-(3.47) of \cite{BGT2}, we obtain
\begin{eqnarray*}
\lefteqn{E\left(\int^{Tt}_0\langle N^{x+h_{k,n}(x)}_r,\varphi\rangle 
dr\right)^3=E
\frac{\partial^3}{\partial\theta^3}v_\theta(x+h_{k,n}(x),Tt)|_{\theta=0}}\\
&=&6E\int^{Tt}_0{\cal T}_{Tt-r}\left(\varphi\int^T_0{\cal T}_{r-u}\left(\varphi\int^u_0{\cal T}_s\varphi ds\right)du\right)(x+h_{k,n}(x))dr\\
&+&3VE\int^{Tt}_0{\cal T}_{Tt-r}\left(\varphi\int^r_0{\cal T}_{r-u}\left(\int^u_0{\cal T}_s\varphi ds\right)^2du\right)(x+h_{k,n}(x))dr\\
&+&6VE\int^{Tt}_0{\cal T}_{Tt-r}\left(\int^r_0{\cal T}_v\varphi dv\int^r_0{\cal T}_{r-u}\left(\varphi\int^u_0{\cal T}_s\varphi ds\right)du\right)(x+h_{k,n}(x))dr\\
&+&3V^2E\int^{Tt}_0{\cal T}_{Tt-r}\left(\int^r_0{\cal T}_v\varphi dv
\int^r_0{\cal T}_{r-u}\left(\int^u_0{\cal T}_s\varphi ds\right)^2du\right)(x+h_{k,n}(x))dr.
\end{eqnarray*}

Without loss of generality we may assume $t=1$.

Repeating the argument as in (\ref{eq:3.5}), (\ref{eq:3.6}), and using the invariance of  Lebesgue measure for ${\cal T}_t$, it is not hard to see that in order to prove (\ref{eq:3.19}) it suffices to show that
\begin{equation}
\label{eq:3.58}
\lim_{T\to\infty}J_i(T)=0,\,\, i=1,2,3,4,
\end{equation}
where
\begin{eqnarray}
\label{eq:3.59}
J_1(T)&=&\frac{1}{F^3_T}\int_{\errita}\int^T_0\phi_2(x)\int^r_0{\cal T}_{r-u}
\left(\phi_2\int^u_0{\cal T}_s\phi_2ds\right)(x)dudrdx,\\
\label{eq:3.60}
J_2(T)&=&\frac{1}{F^3_T}\int_{\errita}\int^T_0\phi_2(x)
\int^r_0{\cal T}_{r-u}\left(\int^u_0{\cal T}_s\phi_2ds\right)^2(x)dudrdx,\\
\label{eq:3.61}
J_3(T)&=&\frac{1}{F^3_T}\int_{\errita}\int^T_0\int^r_0
{\cal T}_v\phi_2(x)dv\int^r_0{\cal T}_{r-u}\left(\phi_2
\int^u_0{\cal T}_s\phi_2ds\right)(x)dudrdx,\\
\label{eq:3.62}
J_4(T)&=&\frac{1}{F^3_T}\int_{\errita}\int^T_0\int^r_0
{\cal T}_v\phi_2(x)dv\int^r_0{\cal T}_{r-u}
\left(\int^u_0{\cal T}_s\phi_2ds\right)^2(x)dudrdx.
\end{eqnarray}

Assume $1/2<\alpha<1$. Denote
$$f(x)=\int^1_0p_r(x)dr,\,\, \tilde{\phi}_T(x)=T^{1/\alpha}\phi_2(T^{-1/\alpha}x).$$
We have
\begin{equation}
\label{eq:3.63}
||\tilde{\phi}_T||_1=||\phi_2||_1,\,\, ||
\tilde{\phi}_T||_2=T^{1/2\alpha}||\phi_2||_2,
\end{equation}
($||\cdot||_p$ denotes the norm in $L^p(\erre)$) and
\begin{equation}
\label{eq:3.64}
||f||_2<\infty,
\end{equation}
since $\alpha>1/2$. The proofs of (\ref{eq:3.58}) use (\ref{eq:3.63}) and 
(\ref{eq:3.64}) together with the Schwarz and Young inequalities. For brevity, we show (\ref{eq:3.58}) for $J_1$ and $J_4$ only.

By (\ref{eq:3.59}) and (\ref{eq:2.12}),
\begin{eqnarray*}
J_1(T)&\leq&T^{-7/2+3/2\alpha}\int_{\errita}\phi_2(x)\int^T_0{\cal T}_u\left(\phi_2\int^T_0{\cal T}_s\phi_2ds\right)(x)dudx\\
&=&T^{-3/2-1/2\alpha}\int_{\errita^3}\tilde{\phi}_T(x)f(x-y)\tilde{\phi}_T(y)f(y-z)\tilde{\phi}_T(z)dzdydx,
\end{eqnarray*}
by obvious substitutions and (\ref{eq:3.1}). Hence
\begin{eqnarray*}
J_1(T)&\leq&T^{-3/2-1/2\alpha}||\tilde{\phi}_T||_2||f*(\tilde{\phi}_T(f*\tilde{\phi}_T))||_2\\
&\leq&T^{-3/2-1/2\alpha}||\tilde{\phi}_T||^2_2||\tilde{\phi}_T||_1||f||^2_2\\
&\leq&C T^{-3/2+1/2\alpha}\to 0,
\end{eqnarray*}
by (\ref{eq:3.63}), (\ref{eq:3.64}), and since $\alpha>1/2$.

Similarly, by (\ref{eq:3.62}) and (\ref{eq:2.12}) we have
\begin{eqnarray*}
J_4(T)&\leq&T^{-7/2+3/2\alpha}\int_{\errita}\int^T_0{\cal T}_v\phi_2(x)dv\int^T_0{\cal T}_u\left(\int^T_0{\cal T}_s\phi_2ds\right)^2(x)dudx\\
&=&T^{1/2-1/2\alpha}||(f*\tilde{\phi}_T)(f*(f*\tilde{\phi}_T)^2)||_1\\
&\leq&T^{1/2-1/2\alpha}||f||^3_2||\phi_2||^3_1\to 0,
\end{eqnarray*}
since $\alpha<1$. The remaining limits in (\ref{eq:3.58}) are obtained similarly.

Now assume $\alpha <1/2$. By (\ref{eq:3.7}) we have
\begin{equation}
\label{eq:3.65}
G\phi_2\in L^\infty(\erre)\cap L^2(\erre)\quad{\rm and}\quad G((G\phi_2)^2)\in L^\infty(\erre)\cap L^2(\erre).
\end{equation}
These properties easily imply (\ref{eq:3.58}). For example,
\begin{equation}
\label{eq:3.66}
J_1(T)\leq\frac{T}{F^3_T}\int_{\errita}\phi_2(x)G(\phi_2G\phi_2)(x)dx\leq\frac{C}{\sqrt{T}},
\end{equation}
and
$$J_4(T)\leq \frac{T}{F^3_T}\int_{\errita}G\phi_2(x)G((G\phi_2)^2)(x)dx\leq\frac{C_1}{\sqrt{T}}.$$

Finally, assume $\alpha=1/2$. In this case we use $G\phi_2\in L^\infty(\erre)$ and
\begin{equation}
\label{eq:3.67}
\frac{1}{\log T}\int_{\errita}\left(\int^T_0{\cal T}_u\varphi(x)du\right)^2dx\leq C\quad{\rm for}\quad T>2,
\end{equation}
see (3.33) in \cite{BGT4}.

It is clear that (\ref{eq:3.66}) also holds (recall that 
$F_T=\sqrt{T\log T}$) and (\ref{eq:3.58})  for $J_2$ and $J_3$ 
follows easily from (\ref{eq:3.67}).

We turn to $J_4$, which requires more work. By (\ref{eq:3.62}), the Schwarz inequality and (\ref{eq:3.67}) we have
\begin{equation}
\label{eq:3.68}
J_4(T)\leq \frac{1}{\sqrt{T}\log T}\sqrt{R(T)},
\end{equation}
where
\begin{eqnarray*}
R(T)&=&\int_{\errita}\left(
\int^T_0{\cal T}_u\left(\int^T_0{\cal T}_s\phi_2ds\right)^2(x)du\right)^2dx\\
&=&\int_{\errita^2}\int^T_0\int^T_0p_{u+u'}(y-z)du'du\left(\int^T_0{\cal T}_s\phi_2(y)ds\right)^2\left(\int^T_0{\cal T}_{s'}\phi_2(2)ds'\right)^2dzdy\\
&=&R_1(T)+R_2(T),
\end{eqnarray*}
where
\begin{eqnarray*}
R_1(T)&=&\int_{\errita^2}\int^1_0\int^T_0\ldots,\\
R_2(T)&=&\int_{\errita^2}\int^T_1\int^T_0\ldots.
\end{eqnarray*}
We have
\begin{eqnarray*}
R_2(T)&\leq&p_1(0)\int^T_1\int^T_0(u+u')^{-2}dudu'\left(\int_{\errita}\left(\int^T_0{\cal T}_s\phi_2(y)ds\right)^2dy\right)^2\\
&\leq&C(\log T)^3,
\end{eqnarray*}
by (\ref{eq:3.67}). Using
$$\int^1_0\int^T_0p_{u+u'}(y-z)dudu'\leq \frac{C}{|y-z|^{1/2}},$$
we obtain
$$R_1(T)\leq R'_1(T)+R''_1(T),$$
where
\begin{eqnarray*}
R'_1(T)&=&C\int\!\!\!\int\limits_{\kern-.5cm |y-z|\leq 1}
\frac{1}{|y-z|^{1/2}}(G\phi_2(y))^2(G\phi_2(z))^2dydz,\\
R''_1(T)&=&C\int\!\!\!\int\limits_{\kern-.5cm |y-z|>1}\left(\int^T_0{\cal T}_s\phi_2(y)ds\right)^2\left(\int^T_0{\cal T}_{s'}\phi_2(z)ds'\right)^2dydz.
\end{eqnarray*}
For $R'_1(T)$ we use (\ref{eq:3.7}) and (\ref{eq:3.5}), obtaining
$$R'_1(T)\leq C_1\int_{\errita}\frac{1}{(1+|y|)^2}dy\int_{|x|\leq 1}\frac{1}{|x|^{1/2}}dx<\infty,$$
and by (\ref{eq:3.67}),
$$R''_1(T)\leq C_2(\log T)^2.$$

Putting these estimates intro (\ref{eq:3.68}) we arrive at
$$J_4(T)\leq C\left(\frac{\log T}{T}\right)^{1/2},$$
which completes the proof of the theorem. $\hfill\Box$

\subsection{Proof of Proposition 2.6}

Since the convergence of finite-dimensional distributions has been proved already, by virtue of the Mitoma theorem \cite{M} it remains to show tightness of $\langle X_T,\varphi\rangle, T>1$, in $C([0,\tau],\erre)$ for any fixed $\varphi\in{\cal S}(\erre), \varphi \geq 0$. To this end we prove
\begin{equation}
\label{eq:3.69}
E(\langle X_T(t),\varphi\rangle-\langle X_T(s),\varphi\rangle)^2\leq C(t-s)^a, s<t\leq \tau,
\end{equation}
for some $a>1$.

By (\ref{eq:3.8}) and (\ref{eq:3.9}) we have
\begin{eqnarray*}
\lefteqn{E(\langle X_T(t),\varphi\rangle-\langle X_T(s),\varphi\rangle)^2}\\
&\leq&\sum_{j\in\zetita}\frac{1}{F^2_T}\int^{j+1}_j\int^{Tt}_{Ts}\int^{Tt}_{Ts}E\langle N^x_r,\varphi\rangle\langle N^x_{r'},\varphi\rangle dr'drdx\\
&=&\frac{1}{F^2_T}\int_{\errita}\int^{Tt}_{Ts}\int^{Tt}_{Ts}E\langle N^x_r,\varphi\rangle\langle N^x_{r'},\varphi\rangle dr'drdx.
\end{eqnarray*}
The latter expression is equal to $E\langle X^{(P)}_T(t)-X^{(P)}_T(s),\varphi\rangle^2$ for $X^{(P)}_T$ corresponding to the system starting from the standard Poisson measure (with or without branching), and we know that (\ref{eq:3.69})
is satisfied for $X^{(P)}_T$ (see subsection 3.1 in 
\cite{BGT1}).

\vglue.5cm
\noindent
{\bf Remark} Recall that Proposition 2.6 refers to a special simple choice of $\nu$. For a general $\nu\in{\cal M}$, the proof of (\ref{eq:3.69}) is similar but slightly more involved. One has to estimate an extra term and use an inequality of the type (\ref{eq:3.6}).

\end{document}